\documentclass[a4paper,oneside]{article}
\usepackage[utf8]{inputenc}
\usepackage[english]{babel}
\pdfoutput=1 
\usepackage[affil-it]{authblk}

\usepackage{amsthm,amsmath,amsbsy,amscd,amssymb,graphicx, epsfig,color,verbatim}

\newcommand{\R}{\mathbb R}
\theoremstyle{plain}
\newtheorem{theorem}{Theorem}

\newtheorem{remark}{Remark}

\definecolor{naranja}{cmyk}{0,0.7,0.8,0}

\begin{document}

\date{}
\small{
\title{Systems of Delay Differential Equations: Analysis of a model with feedback}
\author[1, 2]{Pablo Amster}
\author[1]{Carlos H\'ector Daniel Alliera}
\affil[1]{Departamento de Matem\'{a}tica, FCEyN - 
Universidad de Buenos Aires,
Pabell\'on I, Ciudad Universitaria, Buenos Aires, Argentina,
\textcolor{blue}{pamster@dm.uba.ar calliera@dm.uba.ar, http://cms.dm.uba.ar}}

\affil[2]{IMAS-CONICET}
 \maketitle 
\begin{abstract}

\begin{small}
Using topological degree theory, we prove the existence of positive 
 $T$-periodic solutions of a system of delay differential equations 
for models with feedback arising on regulatory mechanisms in which self-regulation is relevant, e.g. in cell physiology. 
 \end{small} 
\end{abstract}

\section{Introduction}

\begin{figure}[h]
\begin{minipage}{9cm}
\begin{small}
{Self-regulatory models are common in nature, as described e.g. in (\cite{mur}), (\cite{ha}) and (\cite{Gb}).\\
Let us consider a system made up of a number of glands as a motivation. Each gland secretes a hormone that allows secretion in the 
{next} gland, which successively generates another hormone to stimulate the next one and so on. 
In the end, a final hormone is released 
which, by increasing its concentration, will inhibit the secretion of previous hormones that allowed the production process. This generates the decay of this hormone to a minimum threshold that re-activates the cycle again.\\
This behavior can be seen in other biochemical processes, such as enzymatic or bacterial models.\\
Topological degree is a useful tool to find stable 
equilibria in a wide variety of models with constant parameters and, furthermore, 
allows to deduce the existence 
of periodic solutions when the constant parameters 
are replaced by periodic functions.

In this work, we study the existence of solutions for a more general model, namely the following system of delay differential equations:}
\end{small}
\end{minipage}
\ \
\hfill \begin{minipage}{9cm}
\begin{center}
\includegraphics[scale=.7]{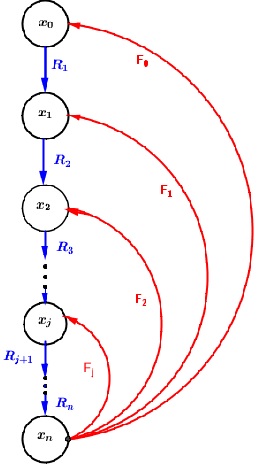} {\caption{A system with feedback}}  
\end{center}
\end{minipage}
\end{figure}

\begin{equation}
\label{t1}
\left\{\begin{array}{l}
\dfrac{dx_0}{dt}= F(t,x_n(t-\tau_0))-b_{0}(x_0(t)), \\
\\
\dfrac{dx_j}{dt}= G_j(t,x_{j-1}(t-\varepsilon_j),x_n(t-\tau_j))-b_{j}(x_j(t)),\quad 1\leq j\leq n-1 \\
\\
\dfrac{dx_n}{dt}= H(t,x_{n-1}(t-\varepsilon_n))-b_{n}(x_n(t))
\end{array}.
\right.
\end{equation}
Here $\tau_i$ and $\varepsilon_j$, with $0\leq i\leq n-1$ and
$1\leq j\leq n$ are fixed positive delays in time.\\

The features of the model read as follows. 

\begin{enumerate}
\item $F,H:\R\times [0,+\infty)\to [0,\infty)$ and $G_j:\R^3\to [0,\infty)$ are continuous and 
$T-$periodic in the first coordinate for some fixed period $T>0$.

\item $b_{i}:[0,+\infty)\to [0,+\infty)$ is a strictly  increasing function with $b_{i}(0)=0$ for $i= 0,\ldots, n$. 

\item $F$ is nonincreasing in its second coordinate and $F(t,x) > 0$ for all $x\ge 0$  and ${\rm Im}(F)\subseteq {\rm Im}(b_0)$. 

\item $H$ is nondecreasing in its second coordinate and $H(t,x)> 0$ for all $x> 0$ and ${\rm Im}(H)\subseteq {\rm Im}(b_n)$. 

\item $G_j$ is nondecreasing in its second coordinate and nonincreasing in its third coordinate, with $G_j(t,x,y)>0$ for $x>0$ and ${\rm Im}(G_j)\subseteq {\rm Im}(b_j)$, for $j=1,\ldots,n-1$. 

\end{enumerate}

We shall prove the existence
of positive $T-$periodic solutions 
for (\ref{t1}), more precisely:

\begin{theorem}
 \label{main}
 Assume that the previous conditions 1-5 hold. Then problem (\ref{t1}) has at least one $T-$periodic solution $u=(x_0,x_1,...,x_n)$ such that 
 $x_k(t)>0$ for all $t$ and all $k$. 
 
\end{theorem}

\section{Existence of positive periodic solutions}

We shall apply the continuation method in the 
the positive cone 
$$\mathcal{K}:=\{u\in C_T: x_0,x_1,...,x_n\ge 0\}$$
of the Banach space of continuous periodic functions
$$C_{T}:=\{u \in C(\mathbb{R},\mathbb{R}^{n+1}):u(t)=u(t+T)\,\hbox{ for all $t$}\},$$
equipped with the standard uniform norm. Consider the linear operator $L:C ^1\cap C_T\to C$ given by $Lu:=u'$ and the nonlinear operator 
$N:\mathcal{K}\to C_T$ 
defined as the right-hand side of system (\ref{t1}).

For convenience, the average of a function $u$ shall be denoted by $\overline u$, namely 
$\overline u:=\frac 1T\int_0^T u(t)\, dt$. 
Also, identifying $\R^{n+1}$ with the subset of constant functions of $C_T$, we may define the function $\phi:[0,+ \infty)^{n+1}\to \R^{n+1}$ given by
$\phi(x) := \overline {Nx}$, that is:

$$
\phi(x_0,x_1,...,x_n)=  \left( \frac 1T\int_0^T F(t,x_n)\, dt - b_0(x_0),
\frac 1T\int_0^T G_1(t,x_1,x_n)\, dt - b_1(x_1),..., \frac 1T\int_0^T H(t,x_{n-1})\, dt - b_{n}(x_n)
\right).
$$

The following continuation theorem can be easily deduced from the standard
topological degree methods (see e.g. (\cite{PA}).

\begin{theorem}

$\label{TC}$
Assume there exists {$\Omega\subset \mathcal{K}^\circ$} open and bounded such that:

\begin{enumerate}
\item[a)] The problem $Lu=\lambda Nu$ has no solutions on $\partial\Omega$ for $0<\lambda<1$.
\item[b)] $\phi(u)\neq 0$ for all {$u\in\partial\Omega \cap\mathbb{R}^{n+1}$}.
\item[c)] {$deg(\phi,\Omega \cap\mathbb{R}^{n+1}, 0)\neq 0$}, where `deg' denotes the Brouwer degree. 

\end{enumerate}
Then $(\ref{t1})$ has at least one solution in {$\overline \Omega$}.
\end{theorem}

In order to apply Theorem \ref{TC} to our problem, let us assume that $u=(x_0,x_1,...,x_n)\in \mathcal{K}$ is a solution of the system $Lu=\lambda Nu$ for 
some $\lambda\in (0,1)$. We shall obtain bounds that 
will yield an appropriate choice of the subset $\Omega$. 

{In the first place, suppose that $x_0$ achieves its absolute maximum $M_0$ at some value $t^*$, then $x_0'(t^*)=0$ and hence
$$
b_{0}(M_0)=F(t^*,x_n(t^*-\tau_0))\le F(t^*,0).
$$
Fixing a constant 
$\mathcal M_0 > \max_{t\in \mathbb R}b_0^{-1}(F(t,0))$,  we conclude that $x_0(t^*)<\mathcal M_0$. 
Next, 
observe that if $x_1$ achieves its absolute maximum $M_1$ at some $t^*$, then 
$$b_{1}(M_1)=G_1(t^*,x_0(t^*-\varepsilon_1),x_n(t^*-\tau_1))
\le G_1(t^*,M_0,0) \le  G_1(t^*,\mathcal M_0,0).$$
Thus, we may fix a constant $\mathcal M_1 >\max_{t\in\mathbb R}
b_1^{-1}(G_1(t,\mathcal M_0,0))$ and hence $M_1<\mathcal M_1$. 
In the same way, for $j=2,\ldots,n-1$ 
we fix constants 
$\mathcal M_j>\max_{t\in\mathbb R}
b_j^{-1}(G_j(t,\mathcal M_{j-1},0))$ so $x_j(t)<\mathcal M_j$ for all $t$. 
For the last equation, we suppose $x_n$ achieves its absolute maximum $M_n$ for some $t^*$, then 
$$
b_{n}(M_n)=H(t^*,x_{n-1}(t^*-\varepsilon_n))\le H(t^*,\mathcal{M}_{n-1}).
$$
Thus we may fix a constant 
$\mathcal M_n > \max_{t\in\mathbb R}b_n^{-1}(H(t,\mathcal{M}_{n-1}))$ and conclude that $x_n(t^*)<\mathcal M_n$.} 

In order to obtain lower bounds, assume firstly that $x_0$ achieves its absolute minimum $m_0$ at some $t_*$, then 
$$
b_{0}(m_0)=F(t_*,x_n(t_*-\tau_0))\ge F(t_*,\mathcal M_n) > 0.
$$
Then, choosing a positive constant 
$\mathfrak m_0 < \min_{t\in\mathbb R} b_0^{-1}(F(t,\mathcal M_0))$ 
it is seen that $m_0>\mathfrak m_0$. 
In the same way, we fix  
$\mathfrak m_j>0$ is such that  
$\mathfrak m_j < b_j^{-1}( G_j(t,\mathfrak m_{j-1}, \mathcal M_n))$ for all $t$ and conclude that 
then $x_j(t)>\mathfrak m_j$ for all $t$ and $1\leq j\leq n-1$. 

Finally, fix a positive constant $\mathfrak m_n$ such that
$\mathfrak m_n < b_n^{-1}(H(t,\mathfrak m_{n-1}))$ for all $t$, then $x_n(t)> \mathfrak m_n$ for all $t$.

In other words, the first 
condition of the continuation theorem is satisfied on 
$$
\Omega:= \{ (x_0,x_1,...,x_n)\in C_T: \mathfrak m_0 < x_0(t) < \mathcal M_0,..., 
\mathfrak m_j < x_j(t) < \mathcal M_j,..., \mathfrak m_n < x_n(t) < \mathcal M_n\, \hbox{ for all $t$,  $1\leq j\leq n-1$}\}.
$$
On the other hand, 
observe that $Q:=\Omega\cap \R^{n+1} =  
(\mathfrak {m}_0,\mathcal{M}_0) \times \ldots 
\times (\mathfrak {m}_n,\mathcal{M}_n)$,
so we shall study the behavior of the mapping $\phi$ over the faces of $Q$. 

Let $x\in Q$ and suppose $x_0 = \mathfrak{m}_0$, then for some $\hat t$
$$
\dfrac{1}{T}\int_{0}^{T}F(t,x_n)\, dt- b_{0}(\mathfrak{m}_0)= 
F(\hat t,x_n) - b_{0}(\mathfrak{m}_0) >  
F(\hat t,x_n) - F(\hat t,\mathcal M_n) \ge 0.
$$
Now suppose $x_0=\mathcal M_0$, then 
$$
\dfrac{1}{T}\int_{0}^{T}F(t,x_n)\, dt- b_{0}(\mathcal{M}_0)= 
F(\tilde t,x_n) - b_{0}(\mathcal{M}_0) <
F(\tilde t,x_n) - F(\tilde t,0)\le 0.
$$
In the same way, for all $j=1,..,n-1$ it is seen that 
$$\dfrac{1}{T}\int_{0}^{T}G_j(t,x_{j-1},x_n)\, dt - b_{j}(\mathfrak{m}_j)= 
G_j(\hat t,x_{j-1},x_n) - b_{j}(\mathfrak{m}_j) \ge
G(\hat t,\mathfrak m_{j-1},\mathcal M_n) - b_{j}(\mathfrak{m}_j) > 0,
$$
$$
\dfrac{1}{T}\int_{0}^{T}G_j(t,x_{j-1},x_n)\, dt - b_{j}(\mathcal{M}_j) = 
G(\tilde t,x_{j-1},x_n) - b_{j}(\mathcal{M}_j) \le
G(\tilde t,\mathcal M_{j-1},0) - b_{j}(\mathcal{M}_j) < 0
$$
and
$$\dfrac{1}{T}\int_{0}^{T}H(t,x_{n-1})\, dt - b_{n}(\mathfrak{m}_n)= 
H(\hat t,x_{n-1}) - b_{n}(\mathfrak{m}_n) \ge
H(\hat t,\mathfrak m_{n-1}) - b_{n}(\mathfrak{m}_n) > 0,
$$
$$
\dfrac{1}{T}\int_{0}^{T}H(t,x_{n-1})\, dt - b_{n}(\mathcal{M}_n) = 
H(\tilde t,x_{n-1}) - b_{n}(\mathcal{M}_n) \le
H(\tilde t,\mathcal M_{n-1}) - b_{n}(\mathcal{M}_n) < 0.
$$
We deduce that the second condition of the continuation theorem is fulfilled. Moreover, 
if we consider the homotopy 
$h:\overline {\mathrm {Q}}\times [0,1]\to \R^{n+1}$ given by
$$h(x,\lambda):=(1-\lambda)(\mathfrak{p}-x)+\lambda \phi(x)$$
where $$
\mathfrak p:=
\left( \dfrac{\mathcal{M}_0+\mathfrak{m}_0}{2}, \ldots , \dfrac{\mathcal{M}_n+\mathfrak{m}_n}{2}\right)
$$
then 
$h \neq 0$ on $\partial \mathrm Q\times [0,1]$. Indeed, if 
$h(x,\lambda)=0$ for some $x = (x_0,x_1,...,x_n)\in \partial\mathrm{Q}$, then we may 
suppose for example that $x_0=\mathcal M_0$ and hence
$$0=h_1(x,\lambda)= (1-\lambda)\dfrac{\overbrace{\mathfrak{m}_0-\mathcal{M}_0}^{<0}}{2}+\lambda\overbrace{\phi_{1}(\mathcal{M}_0,x_1,...,x_n)}^{<0}<0,$$
a contradiction. The other cases follow similarly. 
By the homotopy invariance of the Brouwer degree, we conclude that   
$$deg(\phi,\mathrm{Q},0) = deg({\mathfrak{p} - Id}, \mathrm{Q},0)= (-1)^{n+1}$$
and the proof is complete.

\section{Examples}
\subsection{Model of Testosterone Secretion}
The following system  is
mentioned in \cite{rw} and \cite{DRD}, which is based on a model proposed by Smith \cite{sm} and gave rise to the general model developed in this work.
\begin{figure}[h]
\begin{minipage}{9cm}
\begin{small}
Let us consider the model shown in Figure 2 for the cycle of the Testosterone hormone (see \cite{mur}), where the different variables denote the concentrations at time $t$ of 
the Luteinising Hormone ($LH$),
 which is represented by $R(t)$, from Hypotalamus, Luteinising Hormone 
Releasing Hormone ($LHRH$), represented by $L(t)$, from Pituitary gland and Testosterone Hormone ($TH$) 
from Testes in man, represented by $X(t)$. \\
\\
A general autonomous model describing the biochemical interaction of the hormones $LH$, $LHRH$ and $TH$ in the male is presented.\\
\\
The model structure consists of a negative feedback system of three delay differential equations.
{\begin{remark} It is seen from the model that high levels of $X$ 
affect the concentration of $R$ and $L$.
\end{remark}}
\end{small}
\end{minipage}
\ \
\hfill \begin{minipage}{9cm}
\begin{center}
\includegraphics[scale=.48]{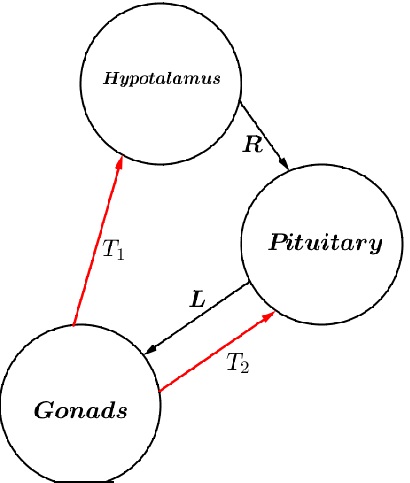}    \caption{Hormone Testosterone cycle}
\end{center}
\end{minipage}
\end{figure}
\begin{equation}
\label{t2}
\begin{array}{l}
\dfrac{dR}{dt}= F(t,X(t-\tau_1))-b_{1}(R(t)), \\
\\

\dfrac{dL}{dt}= g_1 (R(t-\tau_2))-b_{2}(L(t)), \\
\\

\dfrac{dX}{dt}= g_2 (L(t-\tau_3))-b_{3}(X(t))
\end{array}
\end{equation}

{This model has the form of (\ref{t1}) and conditions 1-5 are satisfied if} $b_i (x)$ and $g_j (x)$ 
are positive increasing functions for $i=1, 2, 3$  and $j=1,2$, the delays  
$\tau_i\geq 0$ are constant (at least one of them different from zero) and $F$ is positive, $T$-periodic in $t$ 
and strictly decreasing in $X$.\\
\\
With this structure, Murray \cite{mur} proposed in 1989 a 
simpler (autonomous) system, with:
$$b_i (x)=\beta_i x, \quad \beta_i>0,\quad\quad g_j (x)=\alpha_j x, \quad \alpha_j >0 ,\quad\quad F(x)=\dfrac{\kappa_1}{\kappa_2 + x^{m}},\ m\in\mathbb{N}\quad \tau_1 =\tau_2 =0$$
where $\kappa_j >0$ are constants.\\
The functions $g_j$ represent the rates for productions of $L$ and $X$, $b_i$ are the respective decay rates in the blood stream. It is assumed that each of these hormones is cleared from blood stream according to first order kinetics. (See Das et al (\cite{DRD})



\begin{thebibliography}{1}
\bibitem{PA} {\sc{L. Idels, P. Amster}},
\emph{Existence theorems for some abstract nonlinear non-autonomous systems with delays}, Commun. Nonlinear Sci. Numer. Simulat. 19 (2014) 2974--2982.
\bibitem{DRD} {\sc{P. Das, A. B. Roy and A. Das}},
\emph{Stability and oscillations of a negative feedback delay model for the control of testosterone secretion}, Bio Systems 32 (1994) 61--69.
\bibitem{gk} {\sc{D. Greenhalgh, Q. J. A. Khan}}, {\em A Delay Differential Equation Mathematical Model for the control of the hormonal system of the hypothalamus, the pituitary and the testis in man}. 
Nonlinear Analysis: Theory, Methods and Applications, 71 No. 12 (2009), 925--935. 
\bibitem{mur} {\sc{J. Murray}}, {\em Mathematical Biology. I. An Introduction}, Springer, New York 2001.
\bibitem{rw} {\sc{S. Ruan and J. Wei}}, {\em On the zeros of a third degree exponential polynomial with applications to a delayed model for the control of testosterone secretion}, 
IMA Journal of Math. Applied in Medicine and Biology 18 (2001), 41--52. 
\bibitem{ha} {\sc{S. Hastings, J. Tyson and D. Webster}}, {\em Existence of Periodic Solutions for Negative Feedback Cellular Control Systems}, Journal of Differential Equations 25 (1976), 39--64.
\bibitem{go} {\sc{B. Goodwin}}, \emph{Oscillatory Behaviour in Enzymatic Control Processes}, Adv. Enzyme Regul. 3 (1965), 425--438.
\bibitem{Gb} {\sc{A. Goldbeter}}, \emph{Biochemical Oscillations and Cellular Rhythms}, Cambridge University Press (1996).
\bibitem{sm} {\sc{W. Smith}}, \emph{Hypothalamic regulation of pituitary secretion of luteinizing hormone. II. Feedback control of gonadotropin secretion}, Bulletin of Mathematical Biology, Vol. 42 (1980), 57--78. 
\end{thebibliography}
\end{document}